%

\documentstyle{amsppt}
\loadbold
\magnification=1200
\NoRunningHeads

\topmatter
\title
Finite Combinations of Baire Numbers
\endtitle
\author
Avner Landver
\endauthor
\address
Department of Mathematics,
The University of Kansas,
Lawrence, KS 66045
\endaddress
\email
landver\@kuhub.cc.ukans.edu
\endemail

\abstract
Let $\kappa$ be a regular cardinal.  Consider the
Baire numbers of the spaces $(2^{\theta})_\kappa$ for various 
$\theta \geq \kappa$.  Let $l$ be the number of such different Baire numbers.
Models of set theory with $l=1$ or $l=2$ are known and it is also known that 
$l$ is finite.  We show here that if $\kappa > \omega$, then 
$l$ could be any given finite number. 
\endabstract

\endtopmatter

\document
\baselineskip 18pt
\parskip 6pt

The Baire number of a topological space with no isolated points is the 
minimal cardinality of a family of dense open sets whose intersection 
is empty.  The Baire number (also called the Nov\'ak number [V]) of a partial 
order is the minimal cardinality
of a family of dense sets that has no filter [BS] (i.e. no filter on the given
partial order intersecting all these dense sets non-trivially).
$Fn_\kappa(\theta,2)$ is the collection of all partial functions
$p:\theta \to 2$ such that $|p|<\kappa$, and is partially ordered
by reverse inclusion. For $\kappa$ regular and $\theta \geq \kappa$ we
consider the spaces $(2^{\theta})_\kappa$ whose points are functions from
$\theta$ to $2$ and a typical basic open set is 
$\{f:\theta \to 2 \,\,| \,\, p \subset f \}$ where $p \in Fn_\kappa(\theta,2)$.
We denote the Baire number of $(2^{\theta})_\kappa$ by
${\frak n}_{\kappa}^{\theta}$.  It is not hard to see that 
${\frak n}_{\kappa}^{\theta}$ is also the Baire number of 
$Fn_\kappa(\theta,2)$.  Let us now list some known facts (see [L] \S 1).

\proclaim{Facts}  Let $\kappa$ be a regular cardinal and let 
$\theta \geq \kappa$. Then \newline
1. $\kappa^+ \leq {\frak n}_{\kappa}^{\theta} \leq 2^\kappa$. \newline
2. If $2^{<\kappa} > \kappa$, then ${\frak n}_{\kappa}^{\theta}=\kappa^+$. 
\newline
3. If $\theta_1 \leq \theta_2$, then ${\frak n}_{\kappa}^{\theta_2} \leq
{\frak n}_{\kappa}^{\theta_1}$ and therefore
$\{{\frak n}_{\kappa}^{\theta} : \theta \geq \kappa \,\,
\text {is a cardinal}\}$ is finite. 
\newline
4. If $\theta_1 \leq \theta_2$ and ${\frak n}_{\kappa}^{\theta_2}=\theta_1$,
then ${\frak n}_{\kappa}^{\theta_1}=\theta_1$. \newline
5. If $\theta = {\frak n}_{\kappa}^{2^\kappa}$, then $\theta$ is the unique
cardinal with ${\frak n}_{\kappa}^{\theta}=\theta$ and for every 
$\theta_1 \geq \theta$, ${\frak n}_{\kappa}^{\theta_1}=\theta$.
\endproclaim

A. Miller [M] proved that $cof({\frak n}_{\omega}^{\omega}) > \omega$ but
also produced a model for $cof({\frak n}_{\omega}^{\omega_1}) = \omega$.
In this model $|\{{\frak n}_{\omega}^{\theta} : \theta \geq \omega \,\,
\text {is a cardinal}\}|=2$. Similar models for $\kappa > \omega$ can be
found in [L].  In his above mentioned paper, Miller uses a countable support
product of $Fn_\omega(\omega,2)$ to increase 
${\frak n}_{\omega}^{\omega}$ without
changing the value of ${\frak n}_{\omega}^{\omega_1}$ (and hence getting
${\frak n}_{\omega}^{\omega_1} < {\frak n}_{\omega}^{\omega}$).  
This idea will be used next to prove the following theorem.

\proclaim{Theorem}
Let $\kappa > \omega$ be a regular cardinal.  If ZFC is consistent, then
for every $1 \leq l \in \omega$, ZFC is consistent with 
$|\{{\frak n}_{\kappa}^{\theta} : \theta \geq \kappa \,\,
\text {is a cardinal}\}|=l $.
\endproclaim

This answers ([L] 1.6) for $\kappa > \omega$. 
We do not know whether the Theorem is true for $\kappa = \omega$.
Before we turn to the proof of the theorem, we will need the following lemma
which is due to Miller. The proof of the lemma is a forcing argument
that uses
$\diamondsuit_\kappa$.  The use of $\diamondsuit$'s in forcing arguments
originated in [B]; for other such arguments see [Ka], [L] and [L1].

\definition {Definition}
For the cardinals $\kappa, \theta, \lambda$, let $Q_\kappa(\theta, \lambda)$
be the product of $\lambda$ many copies of $Fn_\kappa(\theta,2)$ with support 
of cardinality $\leq \kappa$.
A condition $q \in Q_\kappa(\theta, \lambda)$ is a function with 
$dom(q) \in [\lambda]^{\leq\kappa}$
and such that for every $\alpha \in dom(q)$, $q(\alpha)\in Fn_\kappa(\theta,2)$.
The partial ordering is defined by putting $q \leq p$ if and only if 
$dom(q) \supset dom(p)$ and for every $\alpha \in dom(p)$, 
$q(\alpha) \supset p(\alpha)$. 
If $\{q_\alpha : \alpha < \gamma \} \subset Q_\kappa(\theta, \lambda)$ 
have a lower bound in $Q_\kappa(\theta, \lambda)$, 
then let us denote the largest lower bound by 
$\bigwedge_{\alpha < \gamma}q_\alpha$.
\enddefinition

\proclaim {Lemma}
Let $\kappa > \omega$ be a regular cardinal such that $\diamondsuit_\kappa$
holds.  Let $\lambda, \theta \geq \kappa$ be cardinals. Let 
$Q = Q_\kappa(\theta, \lambda)$. 
Then forcing with $Q$ over $V$ has the following property:
for every function $f \: \kappa \to V$ in the extension there is a set $A \in V$ 
such that $(|A|=\kappa)^V$ and $range(f) \subset A$ 
(in particular, forcing with $Q$ preserves $\kappa^+$).
\endproclaim

\demo {Proof of the lemma}
Assume that 
$$
q_0\Vdash_Q ``\tau\: \kappa \to V".
$$
Let $M$ be an elementary substructure of the universe such that
$|M|=\kappa$, $M$ is closed under sequences of length $<\kappa$
(i.e. for every $\alpha \in \kappa$, $^\alpha M \subset M$),
and such that $q_0, Q, \lambda , \theta , \kappa , \tau$ are all in $M$.
Notice that every set in $M$ that has cardinality $\leq \kappa$ is also a
subset of $M$. Therefore, if $q \in Q\cap M$, then $q\subset M$.

Let $L=\{\lambda_\xi : \xi < \kappa \} = M\cap \lambda$, and
$T=\{\theta_\xi : \xi < \kappa \} = M\cap \theta$.
For every $\xi < \kappa$, let $L_\xi=\{\lambda_\delta : \delta < \xi \}$, and 
$T_\xi=\{\theta_\delta : \delta < \xi \}$.  Notice that 
$L_\xi, T_\xi \in M$.

For every $\alpha \in \kappa$ we define a function 
$B_\alpha \: Q \to \wp(\alpha \times \alpha)$
as follows:
$$
(\xi, \eta) \in B_\alpha(q) \iff [\lambda_\xi \in dom(q) \land \theta_\eta \in
dom(q(\lambda_\xi)) \land q(\lambda_\xi)(\theta_\eta)=1].
$$
Notice that for every $\alpha \in \kappa$,
$B_\alpha \in M$ (because $L_\alpha, T_\alpha \in M$).

Now, let us fix a $\diamondsuit_\kappa$-sequence 
$I=\{I_\xi : \xi < \kappa \}$ on $\kappa \times \kappa$.  
Notice that for every $\xi < \kappa$, $I_\xi \in M$.
We are now ready to construct a decreasing sequence 
$\{q_\alpha : \alpha < \kappa \} \subset Q\cap M$, 
below $q_0$, that satisfies the following conditions:
\roster
\item $\alpha < \beta \implies q_\beta \leq q_\alpha$.
\item $(\forall \alpha < \kappa)\,\,L_\alpha \subset dom(q_\alpha)$.
\item $\alpha < \beta \implies 
q_\beta \restriction L_\alpha = q_\alpha \restriction L_\alpha$.
\item If $\alpha \in \kappa$ is a limit ordinal, then 
$q_\alpha = \bigwedge_{\beta < \alpha}q_\beta$.  (Notice that 
$q_\alpha \in Q \cap M$ because $\{q_\beta : \beta < \alpha \} \in M$.)
\item Given $q_\alpha$ let us define $q_{\alpha+1}$. \newline
\underbar{Case (i)}: There exist $r \leq q_\alpha$ such that
for every $\xi < \alpha$, $dom(r(\lambda_\xi))=T_\alpha$, and 
$B_\alpha(r) = I_\alpha$, and 
$r$ decides $\tau \restriction \alpha$.
In this case, the same is true in $M$.  Hence there are 
$r_\alpha, t_\alpha \in M$ such that $r_\alpha \leq q_\alpha$,
and for every $\xi < \alpha$, $dom(r_\alpha(\lambda_\xi))=T_\alpha$, and 
$B_\alpha(r_\alpha) = I_\alpha$, and
$$
r_\alpha \Vdash_Q ``\tau \restriction \alpha = t_\alpha".
$$
Let $q_{\alpha+1}$ be defined as follows:   
$q_{\alpha+1} = (q_\alpha \restriction L_\alpha) \cup 
(r_\alpha \restriction (dom(r_\alpha) \setminus L_\alpha))$.\newline
\underbar{Case (ii)}: $\neg$ (case (i)).  Let $q_{\alpha+1} \leq q_\alpha$ be
any extension in $M$ that satisfies (2) and (3),  
and let $t_\alpha = \emptyset$.
\endroster
Finally, define $q = \bigwedge_{\alpha < \kappa}q_\alpha$.
By (1) and (3) of the construction, $q \in Q$.
By (2), $dom(q)=L$.
Let $A = \bigcup_{\alpha < \kappa} range(t_\alpha)$.  We claim that 
$$
q \Vdash_Q ``range(\tau) \subset A".
$$
Assume not.  Let $s \leq q$, and $\delta \in \kappa$ be such that 
$s \Vdash_Q ``\tau(\delta) \notin A"$.
Let us define a decreasing sequence $\{s_\alpha : \alpha < \kappa \}$ 
in $Q$ that satisfies the following conditions:
\roster
\item $s_0=s$.
\item $(\forall \alpha < \kappa)\,\, s_\alpha$ decides $\tau \restriction
\alpha$.
\item If $\alpha < \kappa$ is a limit ordinal, then 
$s_\alpha = \bigwedge_{\beta < \alpha} s_\beta$.
\item $(\forall \alpha < \kappa)(\forall \xi < \kappa)\,\,
T_\alpha \subset dom(s_\alpha(\lambda_\xi))$.
\endroster

Now let $B= \{(\xi, \eta) \in \kappa \times \kappa : 
s_{\eta +1}(\lambda_\xi)(\theta_\eta)=1\}$.
Notice that for every $\alpha < \kappa$, 
$B\cap\alpha\times\alpha = B_\alpha(s_\alpha)$.
Let $C=\{\alpha < \kappa : (\forall \xi < \alpha)\,\,
dom(s_\alpha(\lambda_\xi))=
T_\alpha\}$; $C$ is a club.
In addition, $S=\{\alpha < \kappa : B\cap\alpha\times\alpha=I_\alpha\}$ is
stationary.  Pick $\alpha \in C\cap S$ such that $\alpha > \delta$.
Then $s_\alpha$ witnesses that case (i) of part (5) in the construction of 
$\{q_\alpha : \alpha < \kappa \}$ holds (i.e. $r=s_\alpha$). 
So, we are given 
$r_\alpha, t_\alpha \in M$ such that $r_\alpha \leq q_\alpha$, and
$r_\alpha \Vdash_Q ``\tau \restriction \alpha = t_\alpha"$.  Hence
$$
r_\alpha \Vdash_Q ``\tau(\delta) \in A".
$$
But $s_\alpha \leq r_\alpha$, and $s_\alpha \leq s$, 
and this implies the desired contradiction.
\qed
\enddemo

\demo{Proof of the theorem}
Since the theorem is trivial for $l=1$, let us assume that $l \geq 2$.
Start with a model V of ZFC $+$ GCH $+$ $\diamondsuit_\kappa$. Let
$$
\kappa \leq \theta_1 < \theta_2 < \dots < \theta_l
$$
be cardinals with $\theta_i \neq \kappa^+$, and $\theta_l = \theta_{l-1}^+$,
and such that if $\theta_i \neq \kappa$, then $cof(\theta_i) > \kappa$.
Let
$$
\lambda_1 > \lambda_2 > \dots > \lambda_l = \theta_l
$$
be cardinals with $\lambda_1 = \lambda_2^+$ and such that 
$cof(\lambda_i) > \kappa^+$.

Let $Q_i = Q_\kappa(\theta_i, \lambda_i)$ .
Let us force with 
$$
P = Q_1  \times \dots \times Q_{l-1}.
$$
By the GCH, the partial orders $Fn_{\kappa}(\theta_i,2)$ all have the
$\kappa^+$.c.c.\ ([K]  VII 6.10).  Therefore, 
$P$ is (isomorphic to) a product of $\kappa^+$.c.c.\ partial orders
with support of size $\leq \kappa$.  Now use a delta system lemma and the 
Erd\"os-Rado theorem ($(2^\kappa)^+ \rightarrow (\kappa^+)_{\kappa}^{2}$) 
to show that $P$ is $\kappa^{++}$.c.c.\ ([K] VIII(B7)), and 
hence $P$ preserves cardinals $\geq \kappa^{++}$.  Clearly, 
$P$ is $\kappa$-closed
and therefore cardinals $\leq \kappa$ are preserved.
Finally, by the Lemma, $\kappa^+$ is preserved as well.  

Let $G$ be a $P$-generic filter over $V$.  
Let $\theta \ne \kappa^+$ be a cardinal with 
$\kappa \leq \theta \leq \theta_l $.  
Let $i$ be the minimal such that
$\theta \leq \theta_i$.  Let us show that 
$$
{\frak n}_{\kappa}^{\theta}=\lambda_i. \tag{$*$}
$$
Notice that $(*)$ suffices for the proof of the theorem since it in particular
shows that ${\frak n}_{\kappa}^{\theta_l}=\theta_l$ and therefore by fact 5, 
$(*)$ implies that 
$$
(\forall \theta \geq \theta_l) \,\,{\frak n}_{\kappa}^{\theta}=\lambda_l.
$$
In the remaining case where 
$\theta = \kappa^+$, $(*)$ implies that
${\frak n}_{\kappa}^{\kappa^+} = \lambda_1$ or
${\frak n}_{\kappa}^{\kappa^+} = \lambda_2$.  Therefore, $(*)$ implies that 
$\{{\frak n}_{\kappa}^{\theta} : \theta \geq \kappa \,\,
\text {is a cardinal}\}=\{\lambda_i : 1 \leq i \leq l \} $.
  
Let us first show that ${\frak n}_{\kappa}^{\theta} \geq \lambda_i$. 
By fact 4, we may assume that $1 \leq i < l$.
Notice that since $P$ is $\kappa$-closed, $Fn_\kappa(\theta,2)$ is absolute
and has cardinality $\theta^{<\kappa} \leq \theta_i < \lambda_i$.  
By the product lemma, we may view forcing with
$P$ as forcing with the product 
$\prod \{Q_j: 1 \leq j < l \,\,\text{and}\,\, j \ne i \} \times Q_i$.
Now, by the definition of $Q_i$ and since $\theta \leq \theta_i$, 
it is easy to see that any collection of $< \lambda_i$ many dense subsets
of $Fn_\kappa(\theta,2)$ in $V[G]$, has a filter.

Finally we show that ${\frak n}_{\kappa}^{\theta} \leq \lambda_i$.
Notice that if $i=1$, then this is clear because 
$(2^\kappa = \lambda_1)^{V[G]}$ (to see this use a counting nice names
argument ([K] VII)). So let us assume that $i > 1$ 
and hence $\theta \geq \kappa^{++}$.
In addition we may assume that $\theta$ is regular
(otherwise, if $\theta$ is singular, then it suffices to prove that 
${\frak n}_{\kappa}^{\theta_{i-1}^{++}}\leq \lambda_i$ since
${\frak n}_{\kappa}^{\theta} \leq {\frak n}_{\kappa}^{\theta_{i-1}^{++}}$).

Let us now view forcing with $P$ as forcing 
with $S \times R$, where 
$$
\gather
S = Q_i  \times \dots \times Q_{l-1}\\
\text{and}\\ 
R = Q_1  \times \dots \times Q_{i-1}.
\endgather
$$
Notice that if $i=l$, then $R=P$ and $S$ is the trivial partial order.
Let $H$ be an $S$-generic filter over $V$, and $K$ be an $R$-generic filter
over $V[H]$ such that $V[H \times K] = V[G]$.
For every $a:\theta \to 2$ with $|a| = \kappa$ let us define
$$
D_a = \{t \in Fn_\kappa (\theta ,2): (\exists \xi \in dom(a))\,\, 
t(\xi) \neq a(\xi)\}.
$$
In $V[H]$, define 
$\Cal{D}=\{D_a \,\,|\,\, a\:\theta \to 2 \,\,\text {and}\,\, 
|a| = \kappa \}$.
$\Cal{D}$ is a collection of dense subsets of $Fn_\kappa(\theta , 2)$
and $|\Cal{D}|=\lambda_i$ (because 
$(2^\kappa = \theta^\kappa = \lambda_i)^{V[H]}$).
Let us show that $\Cal{D}$ has no filter in $V[G]$.

Assume, by way of contradiction, that $F \in V[G]$ is a filter for
$\Cal{D}$.  Assume without loss of generality that
$$
\Vdash_{S \times R} ``F\,\,\text{is a filter for}\,\, \Cal{D}".
$$
Let $\tau$ be a $P$-name for $\bigcup F$.  It suffices to find 
$(s,r) \in  S \times R$ and an $S$-name $\pi$ such that
$$
s \Vdash_S ``[\pi:\theta \to 2\,\,\text{and}\,\,
|\pi|=\kappa \,\,\text{and}\,\,
r \Vdash_R ``\pi \subset \tau"]".
$$

We now work in $V$.  For every $\xi \in \theta$, let 
$(s_\xi , r_\xi) \in S \times R$ and $u_\xi \in 2$ be such that 
$$
(s_\xi , r_\xi) \Vdash ``\tau (\xi) = u_\xi".
$$
Consider $\{r_\xi : \xi \in \theta \}$.
Since $\theta \geq \kappa^{++}$ and $\theta$ is regular, we may use the 
delta system lemma to get $X \in [\theta]^\theta$ such that 
$\{dom(r_\xi) : \xi \in X \}$ form a delta system with a root $\Delta$.
Now, since $|Fn_\kappa (\theta_{i-1},2)| = \theta_{i-1} < \theta$ and
$|\Delta| \leq \kappa$, there exists
$Y \in [X]^\theta$ such that 
$\{r_\xi : \xi \in Y \}$ all agree on $\Delta$ 
(i.e.$(\forall \xi, \eta \in Y)\,\,r_\xi \restriction \Delta = 
r_\eta \restriction \Delta$).

Consider $\{s_\xi : \xi \in Y \}$.
Since $S$ is $\kappa^{++}$.c.c.\ there exists $s' \in S$ and a name $\sigma$ 
with
$$
s' \Vdash_S ``\sigma = \{\xi \in Y : s_\xi \in \Gamma \}\,\, and \,\,
|\sigma| = \theta",
$$
where $\Gamma$ is the canonical name for the $S$-generic filter.  
By the Lemma, there exists $A \in [Y]^\kappa$ and $s \leq s'$ such that
$$
s \Vdash_S ``|\sigma \cap A| = \kappa".
$$
Let $\pi$ be an $S$-name for the function whose domain is
$\sigma \cap A$ and such that for every $\xi \in \sigma \cap A$,
$\pi (\xi) = u_\xi$. 
Let $r = \bigcup \{r_\xi : \xi \in A \}$.
Then $r \in R$ (because $A \subset Y$ and $A \in V$), and 
$$
s \Vdash_S ``[\pi:\theta \to 2,\,\,\text{and}\,\,
|\pi|=\kappa, \,\,\text{and}\,\,
r \Vdash_R ``\pi \subset \tau"]". \qed
$$
\enddemo

\remark {Remark 1}
If $\kappa=\omega$, then it is known that $P$ (defined as in the proof of the 
Theorem but for $\kappa =\omega$)
collapses $\omega_1$ ([K] VIII(E4) and [M] p. 280), and (assuming CH) 
is $\aleph_2.c.c$.
What one needs in order to get the argument of the Theorem to go through for 
the case $\kappa =\omega$, is the following:  
if $\sigma$ is a set in the extension that is
unbounded in $(\omega_2)^V$, then there exists a countable set $A$ in $V$
such that $A\cap\sigma$ is infinite.  This is false by the following 
Proposition.
\endremark

\proclaim {Proposition}
Let $\lambda \geq \omega$, and $\theta > \omega$ be cardinals. 
Let $Q = Q_\omega(\theta, \lambda)$. 
Then forcing with $Q$ adds a set $\sigma \subset \theta$, that is unbounded
in $\theta$, and 
such that if $A$ is a countable (in $V$) ground model subset of $\theta$,
then $A \cap \sigma$ is finite.
\endproclaim

\demo {Proof}
For every $n \in \omega$, let $g_n$ be the n'th generic function (i.e.
$g_n \: \theta \to 2$, and $g_n(\alpha)=1$ if and only if there exists
$p$ in the $Q$-generic filter such that $p(n)(\alpha)=1$).
Let $\sigma$ be the set defined in the extension by
$\sigma = \{\alpha \in \theta : (\forall n \in \omega)\,\,g_n(\alpha)=1\}$.
Since $\theta \geq (\omega_1)^V$, and the supports of members of $Q$ are 
countable, it is not hard to see that $\sigma$ is unbounded in $\theta$.
Now let $p \in Q$, and $A \in [\theta]^{\aleph_0}$.
Let us find $q \leq p$ such that
$q \Vdash ``|A \cap \sigma | < \aleph_0"$.
We may assume that $dom(p) \supset \omega$.

Let $A^\ast = \{\alpha \in A : (\exists n \in \omega)\,\,\alpha \notin
dom(p(n)) \}$.  Notice that $A \setminus A^\ast$ is finite.
For every $K \in [\omega]^{<\aleph_0}$ define 
$a(K) =\{\alpha \in A^\ast : (\forall n \notin K)\,\,\alpha \in dom(p(n))\}$;
$a(K)$ is finite.
Fix $\{ \alpha_i : i \in \omega \}$ an enumeration of $A^\ast$.

We now construct $\{q_i : i \in \omega \} \subset Q$, 
$\{n_i : i \in \omega \} \subset \omega$, and
$\{F_i : i \in \omega \}$ finite subsets of $A^\ast$
that satisfy the following conditions:
\roster
\item $q_0 \leq p$ and for every $i \in \omega$, $q_{i+1} \leq q_i$.
\item For every $i \in \omega$, $q_i \restriction (\lambda \setminus 
\{n_k : k \leq i \})
= p \restriction (\lambda \setminus \{n_k : k \leq i \})$.
\item For every $i \in \omega$, $F_i \subset F_{i+1}$, and
$F_i \supset a(\{n_k : k \leq i \})$.
\item $\bigcup_{i \in \omega} F_i = A^\ast$.
\item $i < j \implies q_j \restriction \{n_k : k \leq i \} =
q_i \restriction \{n_k : k \leq i \}$.
\item  For every $i \in \omega$ and every $\alpha \in F_i$,
$q_i \Vdash ``\alpha \notin \sigma"$.
\endroster
\underbar{Stage 0}: Pick $n_0 \in \omega$ with 
$\alpha_0 \notin dom(p(n_0))$.  
Let $F_0 = a(\{n_0\}) \cup \{\alpha_0\}$.
Define $q_0(n_0)$ by: 
$$
q_0(n_0)(\alpha) = \cases 
0 & \alpha \in F_0\\
p(n_0)(\alpha) & \alpha \notin F_0\,\,\text{and}\,\,\alpha \in dom(p(n_0)).
\endcases
$$
\underbar{Stage i+1}: If $\alpha_{i+1} \in F_i$, then $n_{i+1}=n_i, F_{i+1}=F_i$, 
and $q_{i+1}=q_i$.  Otherwise, by (3), 
$\alpha_{i+1} \notin a(\{n_k : k \leq i\})$.
Therefore, we can pick $n_{i+1} \notin \{n_k : k \leq i\}$ such that
$\alpha_{i+1} \notin dom(p(n_{i+1}))$.
By (2), $\alpha_{i+1} \notin dom(q_i(n_{i+1}))$ as well.
Let $F_{i+1} = F_i \cup a(\{n_k : k \leq i+1\}) \cup  \{\alpha_{i+1}\} $.
Define $q_{i+1}(n_{i+1})$ by: 
$$
q_{i+1}(n_{i+1})(\alpha) = \cases 
0 & \alpha \in F_{i+1} \setminus F_i\\
q_i(n_{i+1})(\alpha) &\alpha \notin F_{i+1} \setminus F_i\,\,\text{and}\,\,
\alpha \in dom(q_i(n_{i+1})).
\endcases
$$
Notice that $\alpha \in F_{i+1} \setminus F_i$ implies that either
$\alpha = \alpha_{i+1}$, or $\alpha \in a(\{n_k : k \leq i+1\}) \setminus
a(\{n_k : k \leq i\})$, and in either of these cases 
$\alpha \notin  dom(q_i(n_{i+1}))$.

Finally, let $q = \bigwedge_{i \in \omega}q_i$.
By (2) and (5), $q \in Q$ and clearly, $q \leq p$.
By (4) and (6), 
$q \Vdash ``A^\ast \cap \sigma = \emptyset"$.
\qed
\enddemo

\remark {Remark 2}
In the extension of the above Proposition we also have:
$\sigma$ is an unbounded subset of $\theta$, and
if $ x \in [\sigma]^{\aleph_0}$, then $(\omega_1)^V$ is countable 
in $V[x]$.  This is true because $Q$ is $\aleph_2.c.c.$, and thus
there is $X\in V$ with $|X|=\aleph_1$ and $X \supset x$.  Now one can
enumerate $X$, in $V$, in type $(\omega_1)^V$, and $x$ must be
unbounded in this enumeration since otherwise it would be contained in 
a countable ground model set.

Finally, we would like to mention that the Lemma implies that, the Proposition,
stated for $\kappa > \omega$ (rather than $\omega$), is false.
\endremark

\enddocument